\documentclass[11pt, oneside]{article}

\usepackage{geometry}
\usepackage{booktabs}

\usepackage{amsmath}
\usepackage{amssymb}
\usepackage{algorithm}
\usepackage{algorithmic}

\newcommand{\D}{{\mathrm{d}}}

\newcommand{\bse}{{\boldsymbol{e}}}

\newcommand{\bsx}{{\boldsymbol{x}}}

\newcommand{\bsz}{{\boldsymbol{z}}}

\newcommand{\bsell}{{\boldsymbol{\ell}}}

\newcommand{\bsmu}{{\boldsymbol{\mu}}}

\newcommand{\bsxi}{{\boldsymbol{\xi}}}

\newcommand{\bsrho}{{\boldsymbol{\rho}}}

\newcommand{\bstau}{{\boldsymbol{\tau}}}

\newcommand{\bbE}{{\mathbb{E}}}

\newcommand{\bbR}{{\mathbb{R}}}

\newcommand{\bbV}{{\mathbb{V}}}

\newcommand{\N}{{\mathbb{N}}} 
\newcommand{\R}{{\mathbb{R}}} 

\DeclareSymbolFont{bbold}{U}{bbold}{m}{n}
\DeclareSymbolFontAlphabet{\mathbbold}{bbold}

\newcommand{\calA}{{\mathcal{A}}}

\newcommand{\calI}{{\mathcal{I}}}

\newcommand{\calO}{{\mathcal{O}}}

\usepackage[colorlinks=true,linkcolor=blue,citecolor=blue,urlcolor=black]{hyperref}
\usepackage[font={small,it},labelfont=bf]{caption}

\usepackage{mathtools} 
	
\usepackage{multirow} 
	
\begin{document}

\newcommand{\PRreview}[1]{{\textcolor{black}{#1}}}

\title{A Dimension-Adaptive Multi-Index Monte Carlo Method Applied to a Model of a Heat Exchanger}

\author{
  Pieterjan Robbe\thanks{KU Leuven, Department of Computer Science, Celestijnenlaan 200A bus 2402, 3001 Leuven, Belgium
    (\texttt{\{\href{mailto:pieterjan.robbe@kuleuven.be}{pieterjan.robbe},
    					\href{mailto:dirk.nuyens@kuleuven.be}{dirk.nuyens},
						\href{mailto:stefan.vandewalle@kuleuven.be}{stefan.vandewalle}\}@kuleuven.be}).}
  \and
  Dirk Nuyens\footnotemark[1]
  \and
  Stefan Vandewalle\footnotemark[1]
}

\maketitle

\abstract{We present an adaptive version of the Multi-Index Monte Carlo method, introduced by Haji-Ali, Nobile and Tempone (2016), for simulating PDEs with coefficients that are random fields. A classical technique for sampling from these \PRreview{random fields} is the Karhunen--Lo\`eve expansion. Our adaptive algorithm is based on the adaptive algorithm used in sparse grid cubature as introduced by Gerstner and Griebel (2003), and automatically chooses the number of terms needed in this expansion, as well as the required spatial discretizations of the PDE model. We apply the method to a simplified model of a heat exchanger with random insulator material, where the stochastic characteristics are modeled as a lognormal random field, and we show \PRreview{consistent} computational savings.}

\section{Introduction}\label{PR_sec:Introduction}

\PRreview{A key problem in \emph{uncertainty quantification} is the numerical computation of statistical quantities of interest from solutions to models that involve many random parameters and inputs.} Areas of application include, for example, robust optimization, risk analysis and sensitivity analysis. A particular challenge is solving problems with a high number of uncertainties, leading to the evaluation of high-dimensional integrals. In that case, classical methods such as polynomial chaos~\cite{xiu2009fast,xiu2002modeling} and sparse grids~\cite{bungartz2004sparse} fail, and one must resort to Monte Carlo-like methods.  Recently, an efficient class of such Monte Carlo algorithms was introduced by Giles, see~\cite{barth2011multi,cliffe2011multilevel,giles2008multilevel,giles2015multilevel}. Central to these \emph{multilevel} algorithms is the use of a hierarchy of numerical approximations or \emph{levels}. By redistributing the available computational budget over these levels, taking into account the bias and variance of the different estimators, the error \PRreview{in} the final result is minimized.

A significant extension of the multilevel methodology is the Multi-Index Monte Carlo (MIMC) method, see~\cite{haji2016misc,haji2016multi}. MIMC generalizes the scalar hierarchy of levels to a larger, multidimensional hierarchy of \emph{indices}. This is motivated by the observation that in some applications, changing the level of approximation can be done in several ways, for example in time dependent problems where both time step size and spatial resolution can be varied. Each refinement then corresponds to an index in a multidimensional space. The optimal shape of the hierarchy of indices, based on a priori assumptions on the problem, is analyzed in~\cite{haji2016multi}. However, in most practical problems, such knowledge is not available. \PRreview{Hence the} need for efficient algorithms that automatically detect important dimensions in a problem. Such adaptivity has also been used for deterministic sparse grid cubature in \cite{gerstner2003dimension}. We will develop a similar approach for MIMC.

The paper is organized as follows. In~Sect.~\ref{PR_sec:HeatEquationWithRandomFields}, we introduce a particular example of a PDE with random coefficients: the heat equation with random conductivity. The Multi-Index Monte Carlo method and our adaptive variant are presented in Sect.~\ref{PR_sec:MIMC}. Next, in~Sect.~\ref{PR_sec:HeatExchangerModel}, we introduce a model for a heat exchanger, in which the heat flow is described by the heat equation with random conductivity. We use our adaptive method to compute expected values of the temperature distribution inside the heat exchanger. We show huge computational savings compared to nonadaptive MIMC. We conclude our work in Sect.~\ref{PR_sec:DiscussionAndFutureWork}.

\section{The Heat Equation with Random Conductivity}\label{PR_sec:HeatEquationWithRandomFields}

In this section, we study the linear anisotropic steady state heat equation defined on a domain $D\subset\R^m$, with boundary $\partial D$. The temperature field $T:D\rightarrow\R:\bsx\mapsto T(\bsx)$ satisfies \PRreview{the partial differential equation (PDE)}
\begin{align}\label{PR_eq:DeterministicHeatEquation}
	- \nabla \cdot ( k(\bsx) \nabla T(\bsx)) &= F(\bsx)\quad &\text{for } \bsx \in D,
\end{align}
with $k(\bsx)>0$ the thermal conductivity, $F\in L_2(D)$ a source term, and boundary conditions
\begin{align}
T(\bsx) &= T_1(\bsx) \quad &\text{for } \bsx \in \partial D_1, \nonumber\\
n(\bsx) \cdot (k(\bsx)\nabla T(\bsx)) &= T_2(\bsx) \quad &\text{for } \bsx \in \partial D_2,\nonumber
\end{align}
\PRreview{where $\partial D_1$ and $\partial D_2$ are two disjoint parts of $\partial D$ such that $\partial D = \partial D_1 \cup \partial D_2$. Here, $n(\bsx)$ denotes the exterior unit normal vector to $D$ at $\bsx\in\partial D_2$.}

Consider now the case where equation~\eqref{PR_eq:DeterministicHeatEquation} has a conductivity modeled as a random field, i.e., $k:D\times\Omega\rightarrow\R:(\bsx,\omega)\mapsto k(\bsx,\omega)$ also depends on an event $\omega$ of a probability space $(\Omega, \mathcal{F}, P)$. Then, the solution $T(\bsx,\omega)$ is also a random field and solves almost surely (a.s.)
\begin{align}\label{PR_eq:StochasticHeatEquation}
	- \nabla \cdot ( k(\bsx,\omega) \nabla T(\bsx,\omega)) &= F(\boldsymbol{x}) &\text{for } \bsx \in D \text{ and } \omega\in\Omega,\\
T(\bsx,\omega) &= T_1(\bsx) \quad &\text{for } \bsx \in \partial D_1, \nonumber \\
n(\bsx) \cdot (k(\bsx,\omega)\nabla T(\bsx,\omega)) &= T_2(\bsx) \quad &\text{for } \bsx \in \partial D_2. \nonumber 
\end{align}
\PRreview{For simplicity, we only study the PDE subject to deterministic boundary conditions.}

In what follows, we will develop efficient methods to approximate the expected value
\begin{align}
I(g(\omega))\coloneqq\bbE[g(\omega)] = \int_\Omega g(\omega)\;\D P(\omega),\nonumber
\end{align}
where $g(\omega)=f(T(\cdot,\omega))$ is called the quantity of interest. Typical examples of $g(\omega)$ include the value of the temperature at a certain point, the mean value in (a subdomain of) $D$, or a flux through \PRreview{(a part of) the boundary $\partial D$}.

A commonly used model for the conductivity $k(\bsx,\omega)$ in~\eqref{PR_eq:StochasticHeatEquation} is a lognormal random field, i.e.,
\begin{align}
k(\bsx,\omega)=\exp(Z(\bsx,\omega)),\nonumber
\end{align}
where $Z$ is an underlying Gaussian random field with given mean and covariance. The exponential ensures that the condition $k(\bsx,\omega)>0$ is satisfied for all $\bsx \in D$ and $\omega\in\Omega$, a.s.\ 

In the following, we recall some details about Gaussian random fields that can be found in literature, such as \cite{lord2014introduction,lemaitre2010spectral}. A Gaussian random field $Z(\bsx,\omega)$ is a random field where every vector $\bsz=(Z(\bsx_i,\omega))_{i=1}^M$ follows a multivariate Gaussian distribution with given covariance function for every $\bsx_i\in D$ and $M\in\N$. Specifically, we write $\bsz\sim\mathcal{N}(\bsmu,\Sigma)$, with $\mu_i=\mu(\bsx_i)$ the mean, and with $\Sigma_{i,j}=C(\bsx_i,\bsx_j)\coloneqq\mathrm{cov}(Z(\bsx_i,\omega),Z(\bsx_j,\omega))$ for every $\bsx_i,\bsx_j\in D$, and $C$ the covariance function.

An example of such a covariance function is the \emph{\PRreview{Mat\'ern}} covariance
\begin{align}\label{PR_eq:MaternCovarianceFunction}
C(\bsx_i,\bsx_j) = \sigma^2\frac{1}{2^{\nu-1}\Gamma(\nu)}\left(\sqrt{2\nu}\frac{\|\bsx_i-\bsx_j\|_p}{\lambda}\right)^\nu K_\nu\left(\sqrt{2\nu}\frac{\|\bsx_i-\bsx_j\|_p}{\lambda}\right), \quad \bsx_i,\bsx_j\in D,
\end{align}
where $\Gamma$ is the Gamma function and $K_\nu$ is the modified Bessel function of the second kind. The parameter $\lambda$ is the correlation length, $\sigma^2$ is the (marginal) variance, and $\nu$ is the smoothness of the random field.

Samples of the Gaussian random field can be computed via the \emph{Karhunen--Lo\`eve} (KL) \emph{expansion}
\begin{align}\label{PR_eq:KLExpansion}
	Z(\bsx,\omega)=\mu(\bsx) + \sum_{r=1}^\infty \sqrt{\theta_r}f_r(\bsx) \xi_r(\omega).
\end{align}
In this expansion, the $\xi_r(\omega)$, $r\geq1$, are independent standard normally distributed random numbers and $f_r$ and $\theta_r$ are the solutions to the eigenvalue problem
\begin{align}\label{PR_eq:EigenvalueProblem}
	\int_D C(\bsx_i,\bsx_j)f_r(\bsx_j)\text{d}\bsx_j = \theta_r f_r(\bsx_\PRreview{i}),\quad\bsx_i,\bsx_j\in D\PRreview{,}
\end{align}
\PRreview{where the eigenfunctions $f_r$ need to be normalized for~\eqref{PR_eq:KLExpansion} to hold.} With every event $\omega\in\Omega$ we can associate the (infinite-dimensional) vector $\bsxi(\omega)=\left(\xi_r(\omega)\right)_{r\geq1}$ and, hence, a realization of the random field  $k(\bsx,\omega)$. There exist other methods to generate samples of a random field with given covariance function, such as circulant embedding~\cite{graham2011quasi,lord2014introduction}. Here we choose the KL expansion because of the \emph{best approximation property} described below.

In practice, the infinite sum in~\eqref{PR_eq:KLExpansion} must be truncated after a finite number of terms $s$, that is, $\bsxi(\omega)$ must be truncated to a vector of finite length. The KL expansion gives the best (in MSE sense) $s$-term approximation of the random field if the eigenvalues are ordered in decreasing magnitude~\cite{ghanem1991stochastic,lord2014introduction}. The value of $s$ to reach a certain accuracy depends on the decay rate of the eigenvalues $\theta_r$. The more terms are retained in the expansion, the better the approximation of the random field, but also, the more costly the expansion. This cost involves both the composition of the sum in~\eqref{PR_eq:KLExpansion}, and the (numerical) solution of the eigenvalue problem~\eqref{PR_eq:EigenvalueProblem}. When a lot of terms are required to model the random field, i.e., when the decay of $\theta_r$ is slow, this cost can no longer be ignored compared to the cost of solving the deterministic PDE \PRreview{in every sample of~\eqref{PR_eq:StochasticHeatEquation}}. Hence, it is necessary to construct algorithms that take advantage of the \emph{best approximation property}, and only increase the number of KL terms when required. In~Sect.~\ref{PR_sec:AdaptiveMethod} below, we present an algorithm for such a dimension-adaptive construction of the KL expansion.

\section{The Multi-Index Monte Carlo Method}\label{PR_sec:MIMC}

\PRreview{In sections~\ref{PR_sec:PropertiesOfMonotoneSets} and~\ref{PR_sec:MIMCFormulation} we introduce the Multi-Index Monte Carlo (MIMC) method which was presented and analyzed in~\cite{haji2016multi}. Following that, in section~\ref{PR_sec:AdaptiveMethod} we discuss an adaptive version of the method} based on techniques used in generalized sparse grids, see~\cite{bungartz2004sparse,gerstner2003dimension}. See also~\cite{griebel1992combination} for the \emph{combination technique} on which MIMC is based.

\subsection{Properties of monotone sets}\label{PR_sec:PropertiesOfMonotoneSets}

The formulation of the MIMC method uses the notion of \emph{indices} \PRreview{$\bsell\in S$} and \emph{index sets} $\calI\subseteq S$, where $S\coloneqq\N^d_0=\{\bsell=(\ell_i)_{i=1}^d:\ell_i\in\N_0\}$, with $\N_0=\{0,1,2,\ldots\}$ and $d\geq1$. A \emph{monotone set} is a nonempty set $\calI\subseteq S$ such that for all
\begin{align}\label{PR_eq:DownwardClosed}
\quad\bstau\leq\bsell\in\calI\Rightarrow\bstau\in\calI,
\end{align}
where $\bstau\leq\bsell$ means $\tau_j\leq\ell_j$ for all $j$, see~\cite{chkifa2014high}. Property~\eqref{PR_eq:DownwardClosed} is also known as \emph{downward closedness}. An index set that is monotone is also called a downward closed or \emph{admissible} index set. In the remainder of the text, the index set $\calI$ will always be constructed in such a way that it is an admissible index set. 

Using the definition of the Kronecker sequence \PRreview{$\bse_i\coloneqq(\delta_{ij})_{j=1}^d$}, a monotone set $\calI$ can also be defined using the property
\begin{align}
\big(\bsell\in\calI\quad\mathrm{and}\quad\ell_i\neq0\big)\quad\Rightarrow\quad\bsell-\bse_i\in\calI\quad\mathrm{for\;all\;}\PRreview{i=1,2,\ldots,d}\nonumber.
\end{align}
In other words, for every index $\bsell\ne(0,0,\ldots)$ in a monotone set, all indices with \PRreview{a} smaller \PRreview{(but positive)} \PRreview{entry} in \PRreview{a certain} direction are also included in the set. In the following, we also use the concept of \emph{forward neighbors} of an index $\bsell$, i.e., all indices ${\{\bsell+\bse_i\PRreview{\; : i=1,2,\ldots,d}\}}$, and \emph{backward neighbors} of an index $\bsell$, i.e., all indices $\{\bsell-\bse_i\PRreview{\; : i=1,2,\ldots,d}\}$.

Examples of monotone sets are rectangles
\begin{align}
R(\bsell)
\coloneqq
\left\{\bstau\in S:\bstau\leq \bsell\right\}\nonumber
\end{align}
and simplices
\begin{align}
T_\bsrho(L)
\coloneqq
\left\{\bstau\in S:\bsrho\cdot\bstau\leq L\right\},\nonumber
\end{align}
with $\bsrho\in\R^d_+$ and where $\cdot$ denotes the usual Euclidean scalar product in \PRreview{$\R^d$}.

\subsection{Formulation}\label{PR_sec:MIMCFormulation}

We briefly review the basics of the MIMC method and indicate some of its properties.

Consider the approximation of the expected value of a quantity of interest $g$,
\begin{align}
I(g)
\coloneqq
\bbE[g]
=
\int_\Omega g \;\D P,\nonumber
\end{align}
by an $N$-point Monte Carlo estimator
\begin{align}
Q(g)
\coloneqq
\frac{1}{N}\sum_{n=0}^{N-1} g(\omega_n).\nonumber
\end{align}
Here, the $\omega_n,n=0,1,\ldots$ refer to $N$ random samples from the probability space $\Omega$. \PRreview{Hence}, the estimator \PRreview{itself} is also a random quantity. In our application, the quantity of interest $g$ cannot be evaluated exactly, and we need to resort to discretizations $g_\bsell$, where the different components of $\bsell = (\ell_1, \ldots, \ell_d)$ are different discretization levels of those quantities that need discretization.
Note that the dimensionality of the integral $s$ and the number of discretization dimensions $d$ are not to be confused.

For a given index $\bsell$, define the difference operator in a certain direction $i$, denoted by $\Delta_i$, as
\begin{align}
\Delta_i g_\bsell
\coloneqq
\begin{cases}
g_\bsell - g_{\bsell-\bse_i} &\mbox{if } \ell_i > 0, \\
g_\bsell &\mbox{otherwise},
\end{cases}
\quad i=1,\ldots,d.\nonumber
\end{align}
The MIMC estimator involves a tensor product $\Delta\coloneqq\Delta_1\otimes\cdots\otimes\Delta_d$ of difference operators, where the difference is taken with respect to all backward neighbors of the index $\bsell$.

Using this definition, the MIMC estimator for $I(g)$ can be formulated as
\begin{align}\label{PR_eq:MultiIndexEstimator}
Q_L(g)\coloneqq\sum_{\bsell\in\calI(L)} Q\left(\Delta g_\bsell\right) = 
\sum_{\bsell\in\calI(L)}
\frac{1}{N_\bsell}\sum_{n=0}^{N_\bsell-1} \left(\Delta_1\otimes\cdots\otimes\Delta_d\right) g_\bsell(\omega_{\bsell,n}),
\end{align}
where $\calI(L)$ is an admissible index set. The parameter $L$ governs the size of the index set.

Note that the Multilevel Monte Carlo (MLMC) estimator from~\cite{cliffe2011multilevel,giles2008multilevel,giles2015multilevel} is a special case of the MIMC estimator, where $d=1$. That is, the summation involves a loop over a range of scalar levels $\ell$, and there is no tensor product involved:
\begin{align}
Q_L^\mathrm{(ML)}(g)\coloneqq\sum_{\ell=0}^L Q\left(\Delta g_\ell\right) = 
\sum_{\ell=0}^L
\frac{1}{N_\ell}\sum_{n=0}^{N_\ell-1} \Delta g_\ell(\omega_{\ell,n}).\nonumber
\end{align}

For convenience, we use the following shorthand notation: $E_\bsell\coloneqq|\bbE[\Delta g_\bsell]|$ for the \PRreview{absolute value of the mean} and $V_\bsell\coloneqq\bbV[\Delta g_\bsell]$ for the variance. By $W_\bsell$ we denote the amount of computational work to compute a single realization of the difference $\Delta g_\bsell$. The total work of estimator~\eqref{PR_eq:MultiIndexEstimator} is
\begin{align}\label{PR_eq:TotalWork}
\text{Total Work} = \sum_{\bsell\in\calI(L)} W_\bsell N_\bsell.
\end{align}

In~\eqref{PR_eq:MultiIndexEstimator}, one still has the freedom to choose the index set $\calI(L)$ and the number of samples $N_\bsell$ at each index $\bsell$. In the following, we will show how these two parameters can be quantified.

The objective is to find an index set $\calI(L)$ and sample sizes $N_\bsell$ such that~\eqref{PR_eq:MultiIndexEstimator} achieves a \emph{mean square error} (MSE) smaller than a prescribed tolerance $\epsilon^2$, with the lowest possible cost. From standard statistical analysis, it is known that the MSE can be expressed as a sum of a stochastic error and a discretization error, i.e.,
\begin{align}\label{PR_eq:MeanSquareError}
\bbE\left[\left(Q_L(g)-I(g)\right)^2\right]
&=
\bbE\left[\left(Q_L(g)-\bbE[Q_L(g)]\right)^2\right]
+
\left(\bbE[Q_L(g)]-I(g)\right)^2.
\end{align}
The first term in~\eqref{PR_eq:MeanSquareError} is the variance of the estimator, which, by independence of the events $\omega_{\bsell,n}$, is given by
\begin{align}\label{PR_eq:MIMCVariance}
\bbV[Q_L(g)]=\sum_{\bsell\in\calI(L)}\frac{V_\bsell}{N_\bsell}.
\end{align}
It can be reduced by increasing the number of samples $N_\bsell$. The second term in~\eqref{PR_eq:MeanSquareError} is the square of the bias. It can be reduced by augmenting the index set $\calI(L)$. A sufficient condition to ensure an MSE smaller than $\epsilon^2$, is that both terms in~\eqref{PR_eq:MeanSquareError} are smaller than $\epsilon^2/2$:
\begin{align}
\bbV[Q_L(g)]=\bbE\left[\left(Q_L(g)-\bbE[Q_L(g)]\right)^2\right] &\leq\epsilon^2/2, \quad\text{and} \tag{C1}\label{PR_eq:StatisticalConstraint} \\
|\bbE[Q_L(g)]-I(g)| &\leq\epsilon/\sqrt{2}. \tag{C2}\label{PR_eq:BiasConstraint}
\end{align}
\PRreview{As in~\cite{collier2014continuation} and~\cite{haji2016multi}, we will also use an alternative error splitting, based on a splitting parameter}. The value of this parameter is then computed using a Bayesian approach. This alternative splitting will also be used in our numerical experiments later.

The error splitting in~\eqref{PR_eq:MeanSquareError} will \PRreview{prove} to be essential in the algorithm presented below. Since the total error is the sum of two independent contributions, we can solve for both unknowns $N_\bsell$ and $\calI(L)$ independently. Minimizing the total cost subject to the statistical constraint~\eqref{PR_eq:StatisticalConstraint} \PRreview{will give} the optimal number of samples. Minimizing the total cost subject to the bias constraint~\eqref{PR_eq:BiasConstraint} \PRreview{will yield} the optimal shape of the index set. When using these optimal values for $N_\bsell$, $\calI(L)$, and the error splitting parameter, the cost of the MIMC estimator is minimal\PRreview{,} for a given value of $\epsilon^2$.

\subsubsection{Minimizing the Stochastic Error: Optimal Number of Samples}\label{PR_sec:OptimalNumberOfSamples}

Consider an MIMC estimator with a sufficiently large index set $\calI\PRreview{(L)}$, such that the bias constraint~\eqref{PR_eq:BiasConstraint} is satisfied. Then, one still has to decide the number of samples for each $\bsell\in\calI\PRreview{(L)}$. This freedom can be used to minimize the cost of the MIMC estimator~\eqref{PR_eq:TotalWork} while assuring that the statistical constraint~\eqref{PR_eq:StatisticalConstraint} is satisfied, i.e.,
\begin{align}
&\underset{N_\bsell\in\bbR_+}{\textrm{min}}
\sum_{\bstau\in\calI(L)} N_\bstau W_\bstau \label{PR_eq:OptimizationproblemForOptimalNumberOfSamples} \\
&\textrm{s.t. }
\sum_{\bstau\in\calI(L)}\frac{V_\bstau}{N_\bstau}\leq\frac{\epsilon^2}{2}.\nonumber
\end{align}

This minimization problem can be solved using Lagrange multipliers. The optimal number of samples at each index such that the total cost is minimized, is
\begin{align}\label{PR_eq:OptimalNumberOfSamples}
N_\bsell
=
\frac{2}{\epsilon^2}\sqrt{\frac{V_\bsell}{W_\bsell}}\sum_{\bstau\in\calI(L)}\sqrt{V_\bstau W_\bstau} \quad \textrm{for all }\bsell\in\calI(L).
\end{align}
In practice, this number is rounded up to the nearest integer number of samples. Also, sample variances and estimates for the cost can be used to replace the true variance $V_\bsell$ and true cost $W_\bsell$ at each index. Using~\eqref{PR_eq:OptimalNumberOfSamples}, we can rewrite the total cost of the MIMC estimator as
\begin{align}\label{PR_eq:TotalWork2}
\text{Total Work}=\frac{2}{\epsilon^2}\left(\sum_{\bsell\in\calI(L)}\sqrt{V_\bsell W_\bsell}\right)^2.
\end{align}

\subsubsection{Minimizing the Discretization Error: Optimal Index Sets}

The most simple multi-index method considers indices that are contained in cubes  $\calI(L)=R((L,L\ldots))$ or simplices $\calI(L)=T_{(1,1\ldots)}(L)$. It is possible to extend the latter to the class of general simplices $T_{\bsrho}(\bsell)$. An a priori analysis could then identify important directions in the problem and choose a suitable vector $\bsrho$. However, this approach suffers from two drawbacks. First, such an analysis may be difficult or prohibitively expensive. Furthermore, it is possible that the class of general simplices is inadequate to represent the problem under consideration, especially when mixed directions are involved. In our estimator, we will allow general monotone index sets in the summation~\eqref{PR_eq:MultiIndexEstimator}. \PRreview{The algorithm we designed adaptively} detects important directions in the problem. By a careful construction of the corresponding admissible index set, we hope to achieve an estimator for which the MSE, for a given amount of work, is at least as small as for these classical constructions. Note that as with all adaptive algorithms, the algorithm could be fooled by a quantity of interest for which it seems there is no benefit of extending the index set at some point, and for which essential contributions are hidden at an arbitrary further depth in the index set.

Since the index set is finite, the discretization error is equal to the sum of all neglected contributions, i.e.,
\begin{align}
|\bbE[Q_L(g)]-I(g)|
=
\left|\sum_{\bsell\notin\calI(L)}\bbE[\Delta g_\bsell]\right|
\leq
\sum_{\bsell\notin\calI(L)}E_\bsell.\nonumber
\end{align}
Similar to~\eqref{PR_eq:OptimizationproblemForOptimalNumberOfSamples}, we search for the index set that minimizes the (square root of the) total amount of work~\eqref{PR_eq:TotalWork2}. Here, we impose that the bias constraint~\eqref{PR_eq:BiasConstraint} is satisfied, i.e.,
\begin{align}
&\underset{\calI(L)\subseteq S}{\mathrm{min}}
\sum_{\bsell\in\calI(L)} \sqrt{V_\bsell W_\bsell} \nonumber \\
&\mathrm{s.t.}
\sum_{\bsell\notin\calI(L)}E_\bsell\leq\epsilon/\sqrt{2} \nonumber.
\end{align}

This problem can be formulated as a binary knapsack problem by assigning a \emph{profit indicator} to each index. Define this profit as the ratio of the error contribution and the work contribution, i.e.,
\begin{align}
P_\bsell = \frac{E_\bsell}{\sqrt{V_\bsell W_\bsell}},\label{PR_eq:profitIndicator}
\end{align}
see~\cite{haji2016multi}. A binary knapsack problem is a knapsack problem where the number of copies of each kind of item is either zero or one, i.e., we either include or exclude an index $\bsell$ from the set $\calI(L)$. In the next section, we introduce an adaptive greedy algorithm that solves this knapsack problem, where the profits $P_\bsell$ are used as item weights.

\subsection{An Adaptive Method}\label{PR_sec:AdaptiveMethod}

The goal is to find an admissible index set such that the corresponding MSE is as small as possible \PRreview{subject to an upper bound on the} amount of work. Starting from index $(0,0,\ldots)$, we will successively add indices to the index set such that (a) the resulting index set remains monotone and (b) the error is reduced as much as possible. That is, we require $\calI(0)=\{(0,0,\ldots)\}$ and $\calI(L)\subseteq\calI(L+1)$ for all $L\geq0$. Using the definition of profit above, we can achieve this by always adding the index with the highest profit to the index set. An algorithm that uses this strategy in the context of dimension-adaptive quadrature using sparse grids is presented in~\cite{gerstner2003dimension}. We \PRreview{recall} the main ideas below.

\begin{figure}[htb]
  \centering
  \begin{minipage}{0.875\linewidth}
    \begin{algorithm}[H]
	\caption{\PRreview{Dimension-}Adaptive Multi-Index Monte Carlo}
	\label{PR_alg:AdaptiveMIMC}
	\begin{algorithmic}[]\small
		\STATE{$\bsell\coloneqq(0,\ldots,0)$}
		\STATE{$\calO\coloneqq\varnothing$}
		\STATE{$\calA\coloneqq\{\bsell\}$}
		\STATE{$P_{\bsell}\coloneqq0$}
		\REPEAT
		\STATE{Select index $\PRreview{\bar{\bsell}}$ from $\calA$ with largest profit $P_{\PRreview{\bar{\bsell}}}$}
		\STATE{$\calA\coloneqq\calA\setminus\{\PRreview{\bar{\bsell}}\}$}
		\STATE{$\calO\coloneqq\calO\cup\{\PRreview{\bar{\bsell}}\}$}
		\FOR{$k\;\text{\bf in}\;1,2,\ldots,d$}
		\STATE{$\bstau\coloneqq\PRreview{\bar{\bsell}}+\bse_k$}
		\IF{$\bstau-\bse_j\in\calO$ for all $j=1,2,\ldots,d$ for which $\tau_j>0$}
		\STATE{$\calA\coloneqq\calA\cup\{\bstau\}$}
		\STATE{Take $N^\star$ warm-up samples at index $\bstau$}
		\STATE{Set $Q_\bstau\coloneqq Q(\Delta g_\bstau)$ }
		\STATE{Estimate $V_\bstau$ by~\eqref{PR_eq:MIMCVariance} and $E_\bstau$ by $|Q_\bstau|$}
		\ENDIF
		\ENDFOR
		\FOR{$\bsell\in\calO\cup\calA$}
		\STATE{Compute optimal number of samples $N_\bsell$ using~\eqref{PR_eq:OptimalNumberOfSamples}}
		\STATE{Ensure that at least \PRreview{min$(2,\lceil N_\bsell \rceil)$} samples are taken at each index $\bsell$}
		\STATE{Re-evaluate $Q_\bsell$ and update the estimate of $V_\bsell$ and $E_\bsell$ }
		\STATE{(Re)compute profit indicator $P_\bsell$ using~\eqref{PR_eq:profitIndicator}}
		\ENDFOR
		\UNTIL{$\sum_{\bsell\in\calA} |Q_\bsell|<\epsilon/\sqrt{2}$}
		\RETURN $\sum_{\bsell\in\calO\cup\calA} Q_\bsell$
	\end{algorithmic}
    \end{algorithm}
  \end{minipage}
\end{figure}

The complete algorithm is sketched in~Algorithm~\ref{PR_alg:AdaptiveMIMC}. We assume the current index set $\calI$ is partitioned into two disjoint sets, containing the \emph{active} indices $\calA$ and \emph{old} indices $\calO$, respectively. The active set $\calA$ contains all indices for which none of their forward neighbors are included in the index set $\calI=\calA\cup\calO$. These indices form the boundary of the index set $\calI$ and will actively be adapted in the algorithm. The old index set $\calO$ contains all other indices of the index set, they have at least one forward neighbor in $\calI=\calA\cup\calO$. Equivalently, this means that all backward neighbors of an index in $\calI=\calA\cup\calO$ are always in $\calO$, which means $\calI$ and $\calO$ are admissible index sets. Initially, we set $\calO=\varnothing$ and $\calA=\{(0,0,\ldots)\}$. In every iteration of the adaptive algorithm, the index $\bar{\bsell}$ with the largest profit $P_{\bar{\bsell}}$ is selected from the active set $\calA$. This index is moved from the active set to the old set. Next, all forward neighbors $\bstau$ of $\bar{\bsell}$ are considered. If the neighbor is admissible in the old index set $\calO$, the index is added to the active set $\calA$. A number of warm-up samples are taken at index $\bstau$ to be used in the evaluation of~\eqref{PR_eq:OptimalNumberOfSamples}. After that, we ensure that at least $N_\bsell$ samples are taken at all indices in the index set $\calI=\calA\cup\calO$. Using the updated samples, the profit indicators, as well as the estimates for $V_\bsell$ and $E_\bsell$, are recomputed for all indices in $\calI$. The algorithm continues in the next iteration by selecting the index with the now largest profit, until the condition on the discretization error~\eqref{PR_eq:BiasConstraint} is satisfied. \PRreview{Similar to the approach in~\cite{haji2016multi}}, we use the heuristic bias estimate
\begin{align}
\left|\sum_{\bsell\notin\calI(L)}\bbE[\Delta g_\bsell]\right|
\approx
\sum_{\bsell\in\calA}\left|Q(\Delta g_\bsell)\right|.\label{PR_eq:ComputableBias}
\end{align}
Thus, the absolute value of the Monte Carlo estimators for the differences associated with the indices in the active set $\calA$ act as an estimate for the bias. Finally, note that as soon as an index is added to the active set $\calA$, it is also used in the evaluation of~\eqref{PR_eq:MultiIndexEstimator}. Indeed, it does not make sense to take samples at these indices only to evaluate the profit indicator, and then exclude these samples in the evaluation of the telescoping sum.

\begin{figure}[b]
	\centering
	\hspace{0.85cm}
	\begin{minipage}{0.4\textwidth}
		\caption{Setup for the heat exchanger problem. Hot fluid flows through the left-hand pipe, where a constant heat flux $\Phi_h$ is applied. The cooling fluid in the right-hand pipe has a constant temperature $T_c$. The exterior temperature is $T_e$. The conductivity of the interior conducting material is $k^\mathrm{int}$, while the conductivity of the exterior insulating material is $k^\mathrm{ext}$.}\label{PR_fig:HeatExchanger}
	\end{minipage}
	\hfill
	\begin{minipage}{0.5\textwidth}
		\includegraphics[width=0.8\textwidth]{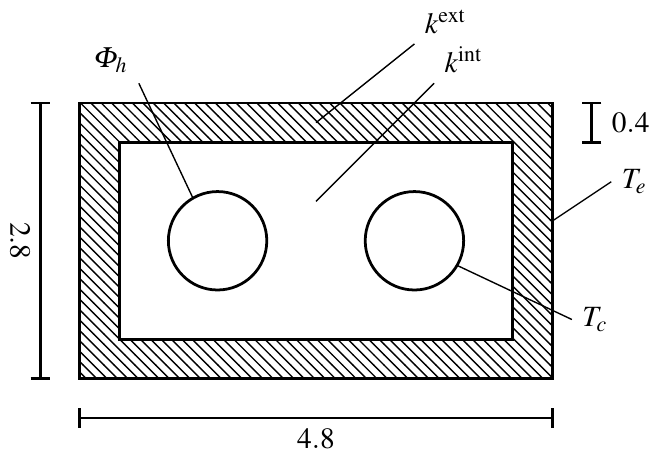}
	\end{minipage}
\end{figure}
\section{A Simple Model for a Heat Exchanger}\label{PR_sec:HeatExchangerModel}

We study the behavior of the adaptive algorithm by applying it to the heat equation with random conductivity from~\eqref{PR_eq:StochasticHeatEquation}. A numerical example, using the strongly simplified model for a heat exchanger from~\cite{lemaitre2010spectral} is presented below. \PRreview{Note that this example, including the choice of its stochastic characteristics, is used for numerical illustration purposes only.}

\begin{figure}[t]
	\centering
	\includegraphics[width=0.44\textwidth]{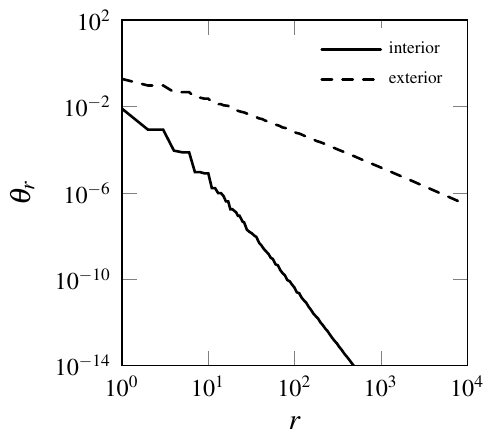}
	\caption{Decay of the eigenvalues $\theta^{\text{int}}$ and $\theta^{\text{ext}}$.}
	\label{PR_fig:EigenvalueDecay}
\end{figure}

\subsection{The Model}\label{PR_sec:TheModel}

We refer to Figure~\ref{PR_fig:HeatExchanger} for a visualization of the description in this section. \PRreview{A} two-dimensional heat exchanger consists of a rectangular piece of material perforated by two circular holes. The first hole contains a hot fluid that injects heat at a constant and known rate $\Phi_h=125/\pi$, and the second hole contains a cooling fluid at a constant temperature $T_c=7.5$. The conductivity of the heat exchanger material is modeled as a lognormal random field $k^\mathrm{int}=\exp(Z^\mathrm{int})$, where $Z^\mathrm{int}$ is a Gaussian random field with mean $\mu^\mathrm{int}=0$ and Mat\'ern covariance with correlation length $\lambda^\mathrm{int}=1$, standard deviation $\sigma^\mathrm{int}=\sqrt{0.1}$, norm $p=1$ and smoothness $\nu^\mathrm{int}=1$, see~\eqref{PR_eq:MaternCovarianceFunction}.

A layer of insulator material is added to the heat exchanger, to thermally insulate it from its surroundings, which has a constant temperature $T_e=20$. The conductivity of the insulator material is modeled as a lognormal random field $k^\mathrm{ext}=\exp(Z^\mathrm{ext})$, where $Z^\mathrm{ext}$ is a Gaussian random field with mean $\mu^\mathrm{ext}=\log(0.01)$ and Mat\'ern covariance with correlation length $\lambda^\mathrm{ext}=0.3$, standard deviation $\sigma^\mathrm{ext}=1$, norm $p=1$ and smoothness $\nu^\mathrm{ext}=0.5$.

\begin{figure}[h!]
	\centering
	\begin{tabular}{cccl}
		\includegraphics[width=0.3\textwidth]{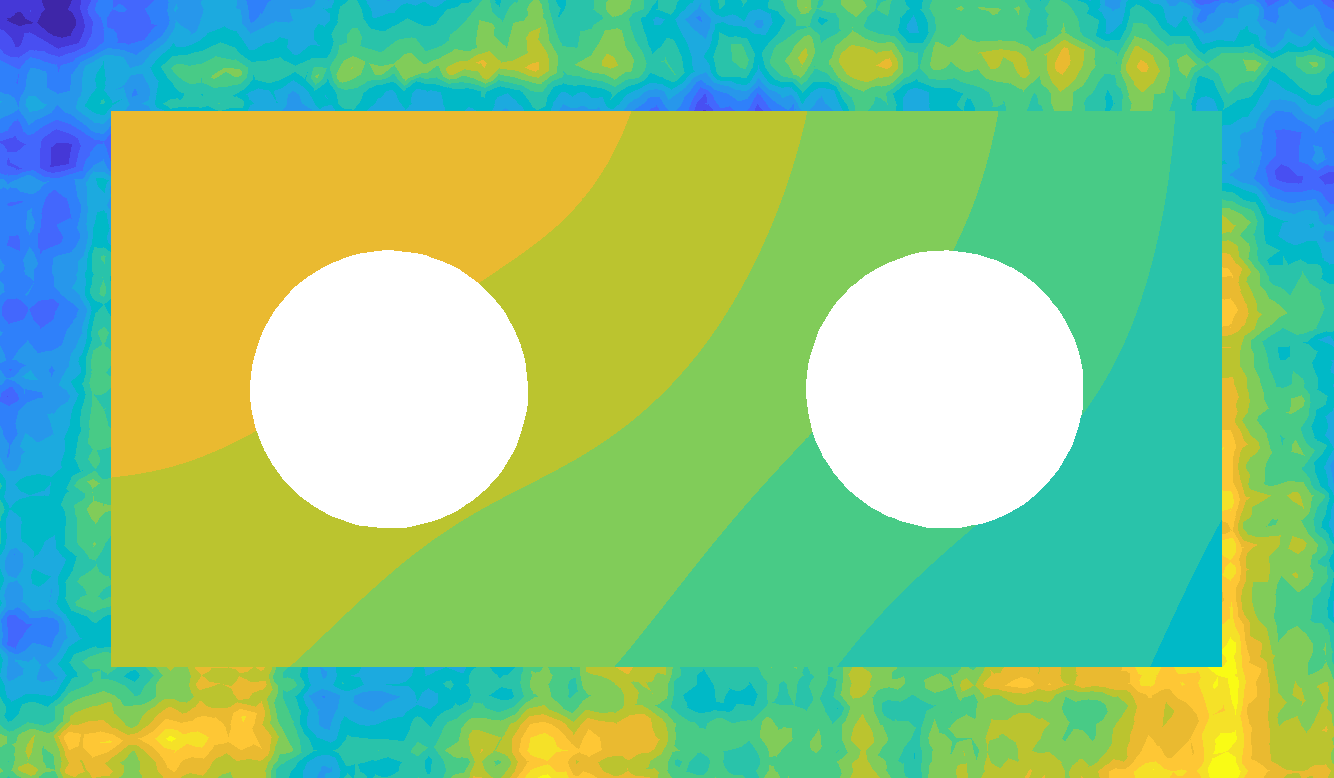} &  
		\includegraphics[width=0.3\textwidth]{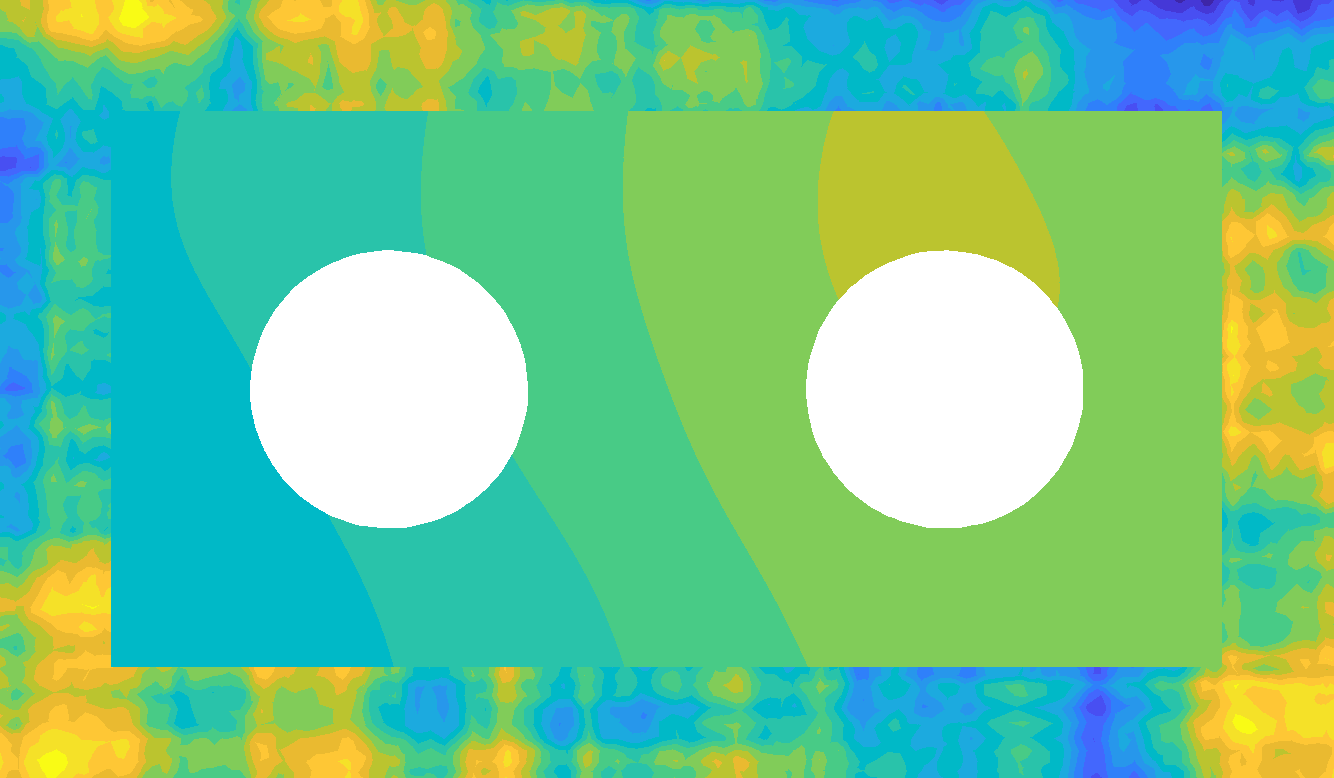} &  
		\includegraphics[width=0.3\textwidth]{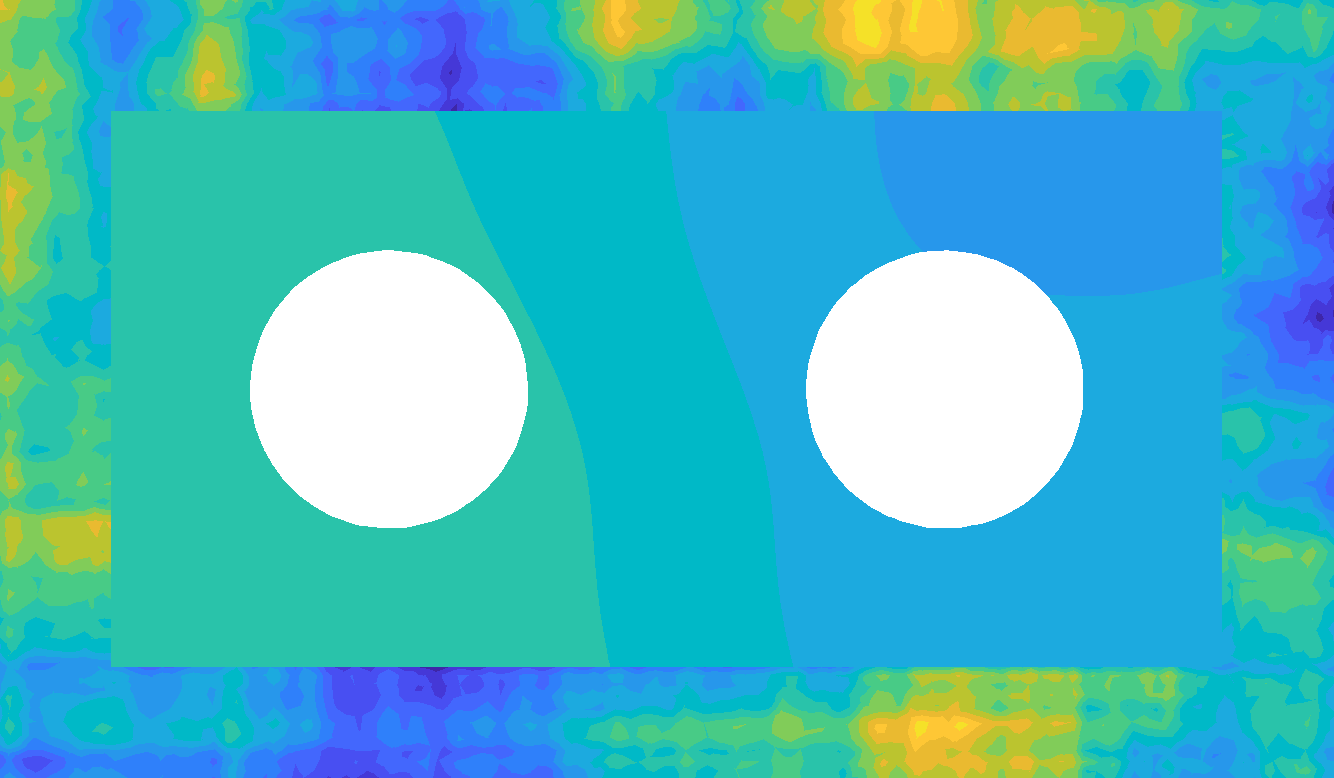} &
		\multirow{1}{*}[0.1665\textwidth]{
			\hspace{-0.45cm}
			\includegraphics[width=0.0795\textwidth]{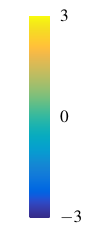}
		}
	\end{tabular}
	\caption{Example realizations of the (zero-mean) Gaussian fields \PRreview{$Z^{\text{int}}$ and $Z^{\text{ext}}$} used in the heat exchanger problem. The insulator material has a lower correlation length and smoothness, and a higher variance. The associated conductivity is \PRreview{$k=\exp(Z^{\text{int}}\cup Z^{\text{ext}})$}. The number of terms used in the KL expansions for $Z^{\text{int}}$ and $Z^{\text{ext}}$ is \PRreview{$512$} and \PRreview{$8\,192$}, respectively.}
	\label{PR_fig:ExampleFields}
\end{figure}

Samples of both random fields are generated using a truncated KL expansion, see~\eqref{PR_eq:KLExpansion}. Figure~\ref{PR_fig:EigenvalueDecay} shows the decay of the two-dimensional eigenvalues for both the conductor (interior) and insulator (exterior) material. These eigenvalues and corresponding eigenfunctions are computed once for the maximal number of terms allowed in the expansion. Every realization of the conductivity $k=\exp(Z)$ is formed using a sample of the Gaussian random fields  $Z^\mathrm{int}$ and $Z^\mathrm{ext}$. Three samples of the (Gaussian) random field $Z$ are shown in~Figure~\ref{PR_fig:ExampleFields}. Note that $Z^\mathrm{int}$ only varies mildly in comparison to $Z^\mathrm{ext}$.

\begin{figure}[h!]
	\centering
	\begin{tabular}{ccc@{\hskip 0.375cm}}
		\includegraphics[width=0.3\textwidth]{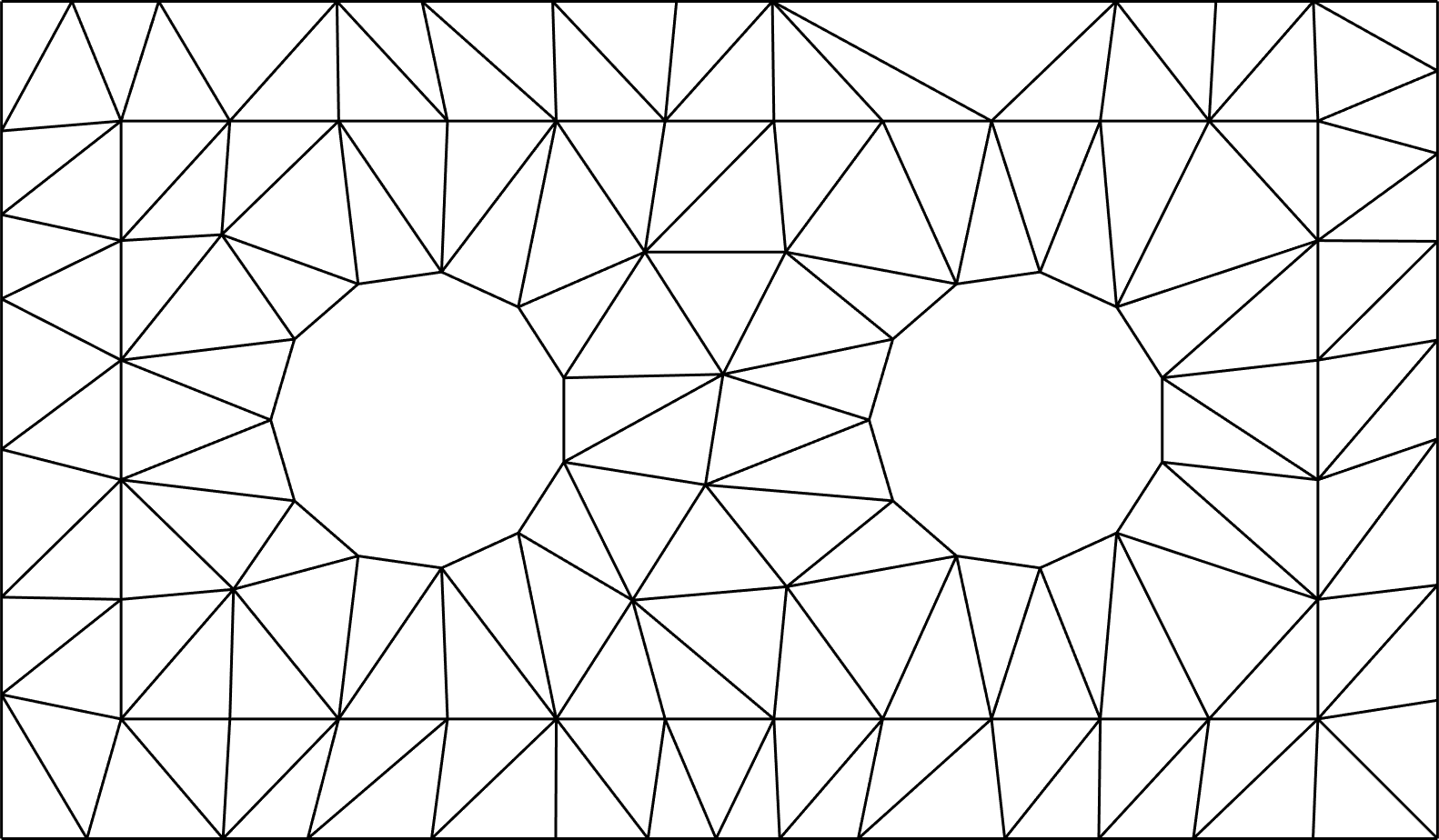} &
		\includegraphics[width=0.3\textwidth]{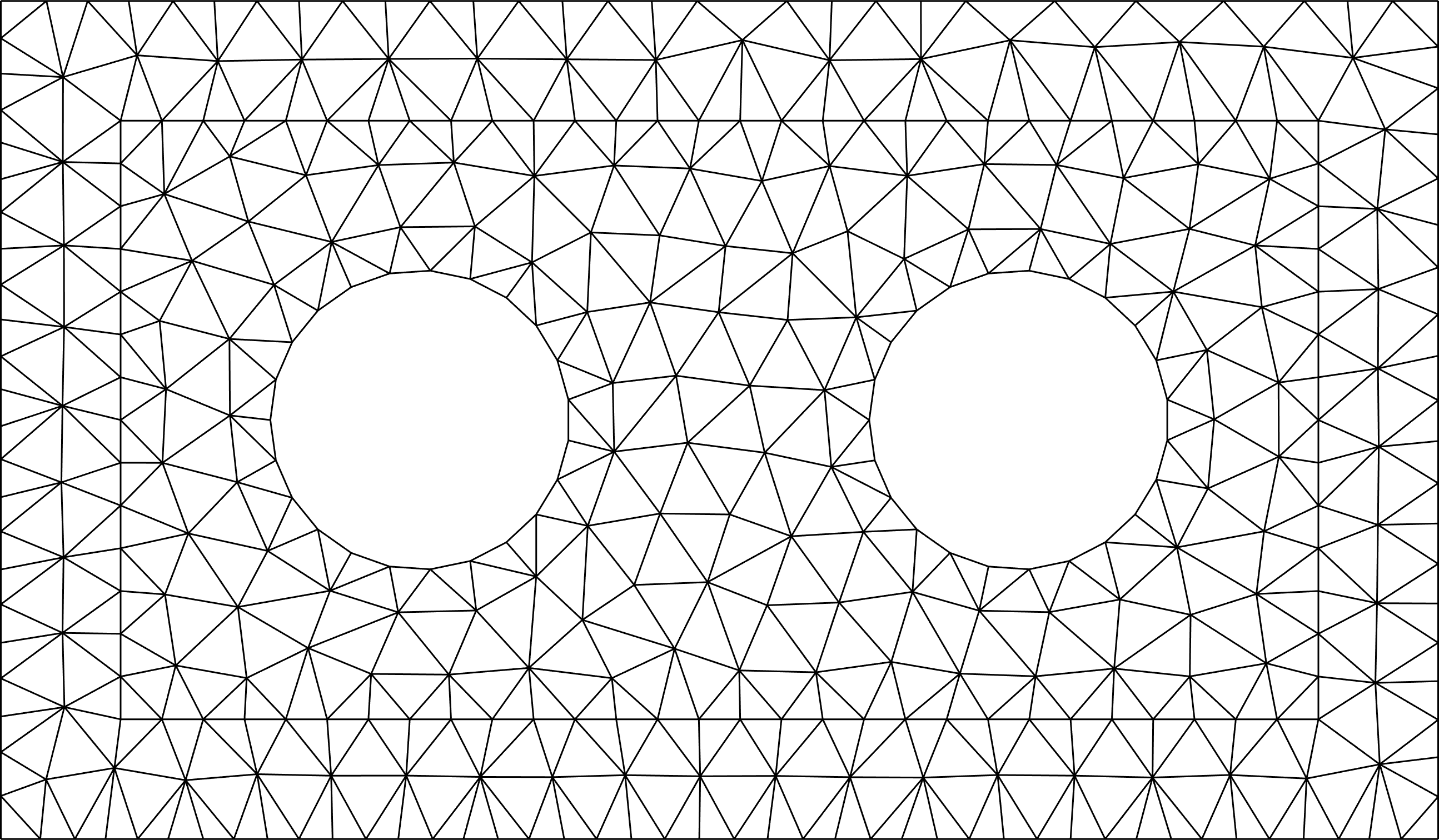} &  
		\includegraphics[width=0.3\textwidth]{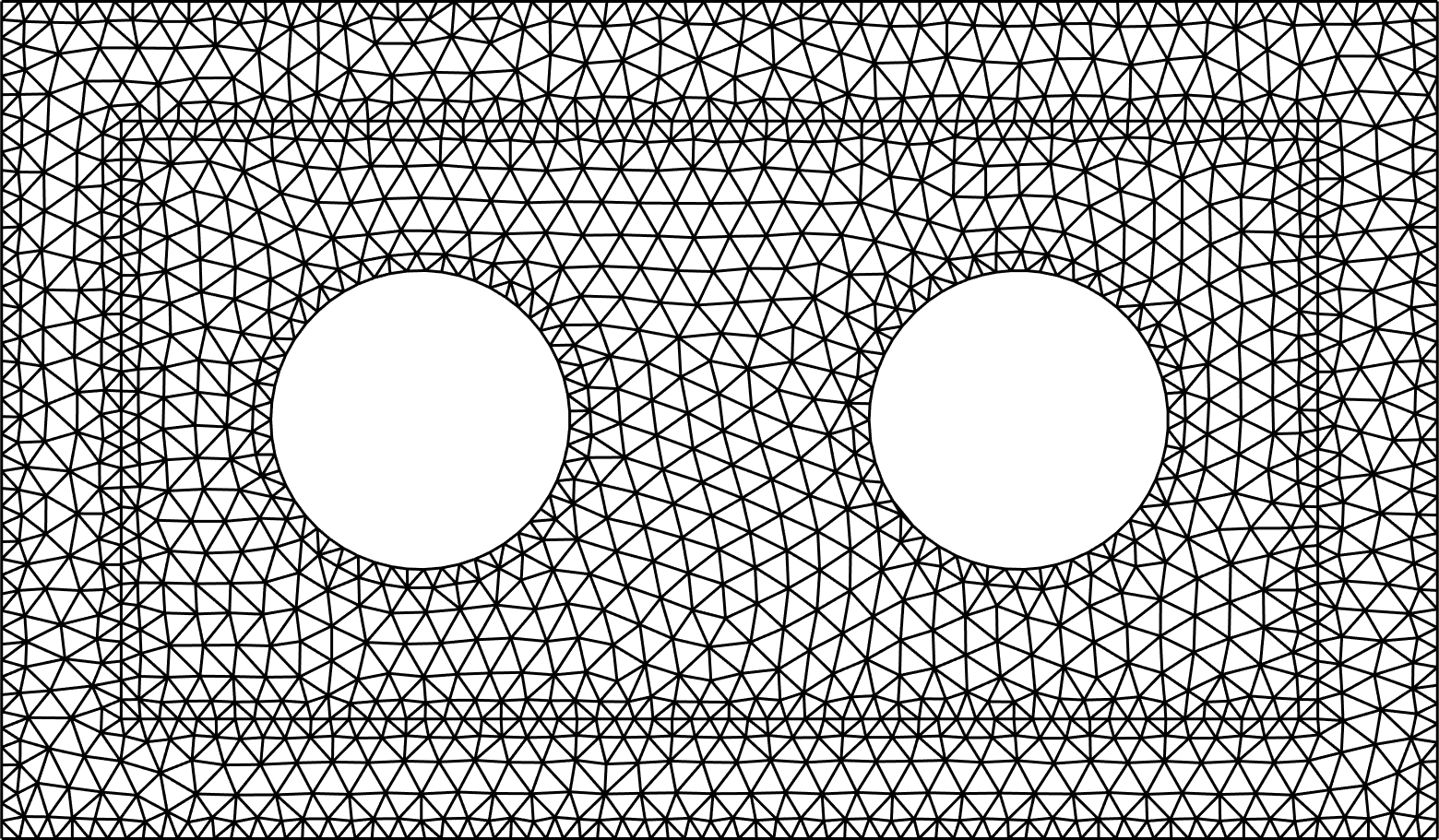} \\
		the coarsest mesh & an intermediate mesh & a fine mesh
	\end{tabular}
	\caption{Some finite-element meshes used in the heat exchanger problem. The coarsest mesh has 102 points (144 elements), the intermediate mesh has 309 points (640 elements), and the fine mesh has 1\,619 points (2\,887 elements). The finest mesh used in the simulations is not shown.}
	\label{PR_fig:FEMeshes}
\end{figure}

For the spatial discretization, we use eleven different \PRreview{nonnested} finite-element (FE) meshes with an increasing number of elements. For every mesh, the number of points is roughly twice the number of points of its predecessor. That way, the size of the finite-element system matrix doubles between successive approximations. The coarsest mesh has 102 points (144 elements), and the finest mesh has 94\,614 points (186\,268 elements). Three examples are shown in~Figure~\ref{PR_fig:FEMeshes}.

The heat flow through the exchanger is described by~\eqref{PR_eq:StochasticHeatEquation}, with source term $F\coloneqq0$. As a quantity of interest, we consider the value of the temperature \PRreview{at} the leftmost point on the boundary of the hot fluid pipe. As shown in~Figure~\ref{PR_fig:MeanField}, this corresponds to the highest expected temperature in the heat exchanger. Note that we have made sure that this point is included on every FE mesh, to avoid an interpolation error.

\subsection{Numerical Results}\label{PR_sec:NumericalResults}

We set up an adaptive MIMC algorithm with three refinement dimensions, i.e. $d=3$. The first dimension corresponds to the spatial discretization, the second dimension is the number of terms in the KL expansion of the conductor material, and the last dimension is used for the number of terms in the KL expansion of the insulator material. The number of terms in \PRreview{either KL expansion doubles} between subsequent approximations, similar to the connection between the different spatial discretizations. If the effect of adding more KL terms to the approximation of the quantity of interest was known \PRreview{in advance}, one could derive the optimal relation between the different approximations, similar to~\cite{haji2014optimization}. This relation will, \PRreview{amongst others}, depend on the decay rate of the eigenvalues of the KL expansion, hence, it will be different for the insulator and conductor material. However, in the absence of this knowledge, doubling the number of terms (a geometric relation, following~\cite{haji2014optimization}) seems \PRreview{an} obvious thing to do. Note that

\noindent the slow eigenvalue decay rate for the insulator material in~Figure~\ref{PR_fig:EigenvalueDecay} is reflected in the number of terms used in the coarsest approximation: index $(\cdot,0,0)$ corresponds to an approximation using \PRreview{$s_0^\mathrm{int}=4$} terms in the KL expansion of the heat exchanger material and \PRreview{$s_0^\mathrm{ext}=64$} terms in the expansion of the insulator material.

\begin{figure}[p]
\includegraphics[width=0.1\textwidth]{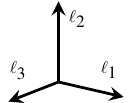}
\hspace{2em}
\scalebox{0.35}{\def\boxcolor{white!90!black}\includegraphics{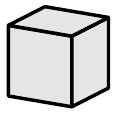}} \raisebox{0.36\baselineskip}{old set}\hfill
\scalebox{0.35}{\def\boxcolor{orange!50!white}\includegraphics{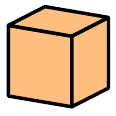}} \raisebox{0.36\baselineskip}{active set}\hfill
\scalebox{0.35}{\def\boxcolor{blue!50!white}\includegraphics{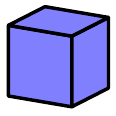}} \raisebox{0.36\baselineskip}{highest profit\hspace{0.4cm}}
 \vspace{\baselineskip} \\
\centering
\begin{tabular}{ccc}\noalign{\smallskip}\hline\noalign{\smallskip}
$L=3$ & $L=4$ & $L=5$ \\ \noalign{\smallskip}\hline\noalign{\smallskip}\noalign{\smallskip}
\includegraphics[width=3cm]{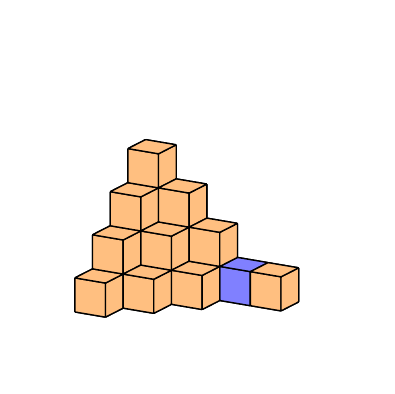}\hphantom{---}&
\includegraphics[width=3cm]{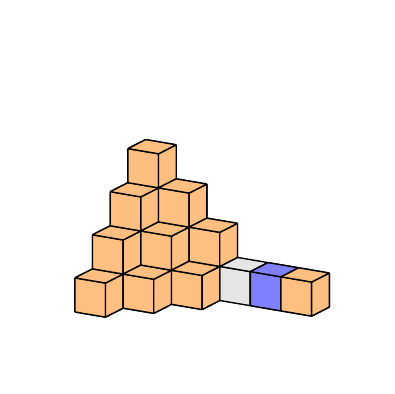}\hphantom{---}&
\includegraphics[width=3cm]{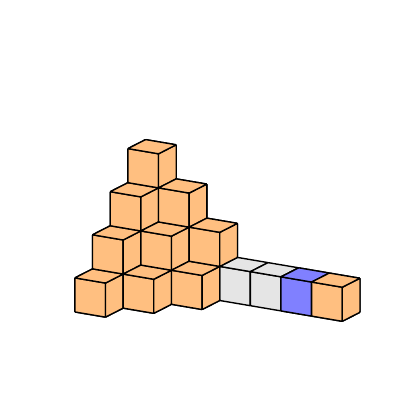}\hphantom{---}\\[1ex]
\end{tabular}
\begin{tabular}{ccc}\noalign{\smallskip}\hline\noalign{\smallskip}
$L=8$ & $L=12$ & $L=16$ \\ \noalign{\smallskip}\hline\noalign{\smallskip}\noalign{\smallskip}
\includegraphics[width=3cm]{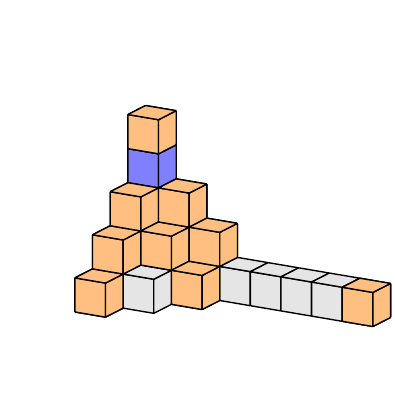}\hphantom{---}&
\includegraphics[width=3cm]{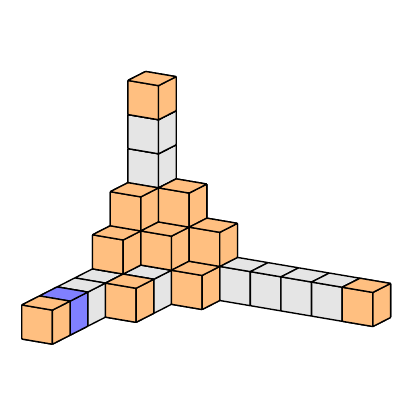}\hphantom{---}&
\includegraphics[width=3cm]{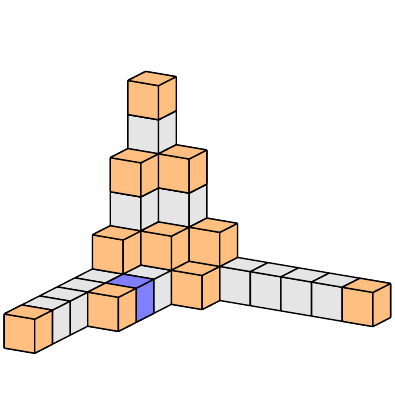}\hphantom{---}\\[1ex]
\end{tabular}
\begin{tabular}{ccc}\noalign{\smallskip}\hline\noalign{\smallskip}
$L=23$ & $L=25$ & $L=29$ \\ \noalign{\smallskip}\hline\noalign{\smallskip}\noalign{\smallskip}
\includegraphics[width=3cm]{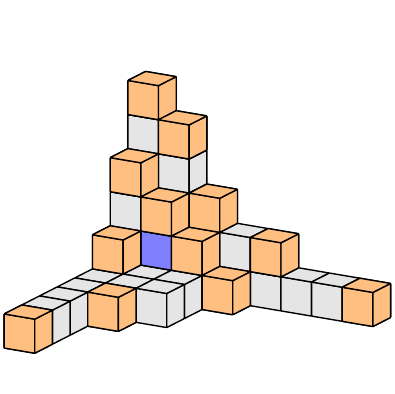}\hphantom{---}&
\includegraphics[width=3cm]{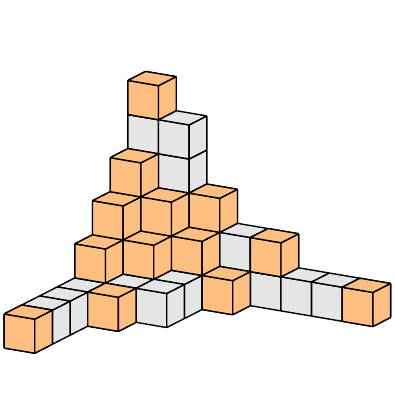}\hphantom{---}&
\includegraphics[width=3cm]{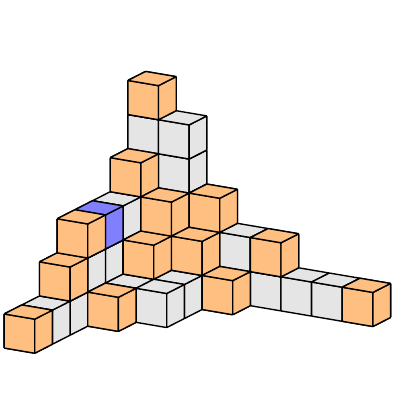}\hphantom{---}\\[1ex]
\end{tabular}
\begin{tabular}{ccc}\noalign{\smallskip}\hline\noalign{\smallskip}
$L=33$ & $L=37$ & $L=40$ \\ \noalign{\smallskip}\hline\noalign{\smallskip}\noalign{\smallskip}
\includegraphics[width=3cm]{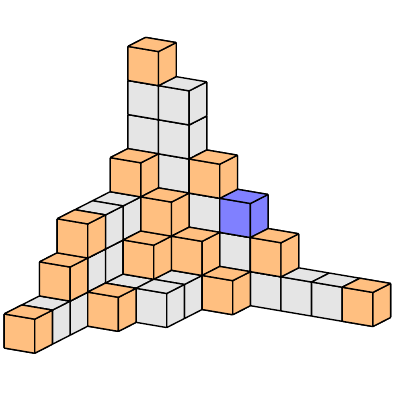}\hphantom{---}&
\includegraphics[width=3cm]{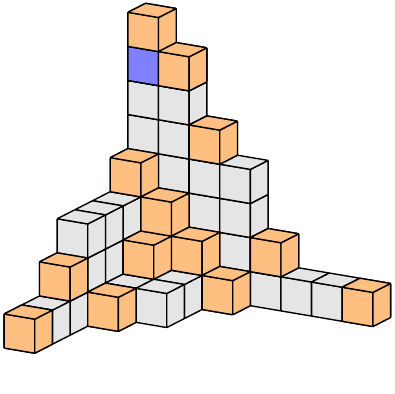}\hphantom{---}&
\includegraphics[width=3cm]{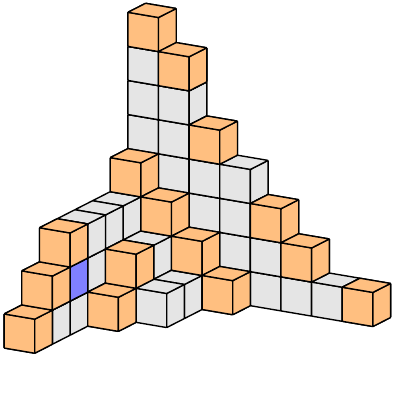}\hphantom{---}\\[1ex]
\end{tabular}
\caption{\PRreview{Examples of nontrivial index sets in the heat exchanger problem for selected iterations in the adaptive algorithm.}}	
\label{PR_fig:IndexSets}
\end{figure}

In practice, we do not start the algorithm from index $(0,\ldots,0)$ as is indicated in~Algorithm~\ref{PR_alg:AdaptiveMIMC}, but start with an index set $T_{(1,1,1)}(2)$, to ensure the availability of robust estimates for the profit indicator on the coarsest approximations. 

\PRreview{The total cost of the computation of $G_\bsell$ is equal to the sum of the cost of composing the random field using the KL expansion and the cost of the finite-element computation. For a given index $\bsell=(\ell_1,\ell_2,\ell_3)$, we assume that there are \texttt{elements}$(\ell_1)$ elements and \texttt{nodes}$(\ell_1)$ nodes in the discretization. The KL expansions at that index use $s_0^\mathrm{int}2^{\ell_2}$ terms for the conductor and $s_0^\mathrm{ext}2^{\ell_3}$ terms for the insulator. We propose the cost model
\begin{align}
C_1(\text{\texttt{elements}}(\ell_1))(s_0^\mathrm{int}2^{\ell_2}+s_0^\mathrm{ext}2^{\ell_3})+C_2(\text{\texttt{nodes}}(\ell_1))^\gamma,
\end{align}
for some suitable constants $C_1$, $C_2$ and $\gamma$. We numerically found the values $C_1=1.596\mathrm{e}{-8}$, $C_2=1.426\mathrm{e}{-6}$ and $\gamma=1.664$. There is no cost involved in computing the quantity of interest $G_\ell$ from the solution $T(\bsx,\cdot)$, since no interpolation is required. The cost $W_\bsell$ of computing a single sample of $\Delta G_\bsell$ can be computed by expansion of the tensor product $\Delta=\Delta_1\otimes\Delta_2\otimes\Delta_3$.} Note that it is also possible to use actual simulation times as measures for the cost. However, since the cost estimate appears in the profit indicator, and thus determines the shape of the index set, one should ensure that the estimates are stable and reliable. Finally, we use the continuation approach from~\cite{collier2014continuation} and run the MIMC algorithm for a sequence of larger tolerances than required, to obtain more accurate estimates of the sample variances, and to avoid having to take warm-up samples at every index.

\begin{table}[t]
	\centering
	\scriptsize
	\caption{\PRreview{Mean, RMSE and running time for MIMC on the left, and for adaptive MIMC on the right. The nonadaptive version uses simplices $T_{(1,1,1)}(L)$ as index set.}}\label{PR_tab:RunTimeComparison}
	\begin{tabular}{p{1.2cm}ccrlccrl}
	\toprule
		& \multicolumn{3}{c}{nonadaptive MIMC} & & \multicolumn{3}{c}{adaptive MIMC} & \\
		$\epsilon_\mathrm{rel}$ & \makebox[1.4cm][c]{mean} & \makebox[1.3cm][c]{\text{RMSE}} & \makebox[1.3cm][r]{time (s)} & \makebox[0.0cm][c]{} & \makebox[1.4cm][c]{mean} & \makebox[1.4cm][c]{RMSE} & \makebox[1.3cm][r]{time (s)} & \makebox[0.0cm][c]{}\\
		\midrule
 2.890$\mathrm{e}{-2}$	& 135.20	& 1.948$\mathrm{e}{0\phantom{-}}$	& 622			& & 133.77	& 2.107$\mathrm{e}{0\phantom{-}}$	& 591			&\\
 1.927$\mathrm{e}{-2}$	& 134.46	& 1.293$\mathrm{e}{0\phantom{-}}$	& 1\,175		& &  132.56	& 1.545$\mathrm{e}{0\phantom{-}}$	& 1\,228			&\\
 1.285$\mathrm{e}{-2}$	& 134.10	& 1.100$\mathrm{e}{0\phantom{-}}$	& 2\,667		& &  132.47	& 1.412$\mathrm{e}{0\phantom{-}}$	& 1\,228			&\\
 8.564$\mathrm{e}{-3}$	& 133.04	& 9.767$\mathrm{e}{-1}	$& 11\,951		& &  132.38	& 1.115$\mathrm{e}{0\phantom{-}}$	& 7\,034			&\\
 5.710$\mathrm{e}{-3}$	& 133.76	& 4.430$\mathrm{e}{-1}$	& 30\,552		& &  133.20	& 4.744$\mathrm{e}{-1}$	& 20\,725			&\\
 3.806$\mathrm{e}{-3}$	& 133.74	& 4.321$\mathrm{e}{-1}	$& 38\,997		& &  133.61	& 1.027$\mathrm{e}{-1}$	& 28\,335			&\\
 2.538$\mathrm{e}{-3}$	& 133.72	& 3.521$\mathrm{e}{-1}	$& 92\,223		& &  133.63	& 1.723$\mathrm{e}{-1}$	& 81\,711			&\\
 1.692$\mathrm{e}{-3}$	& 133.73	& 2.789$\mathrm{e}{-1}$	& 257\,698	& &  133.70	& 1.396$\mathrm{e}{-1}$	& 233\,458			&\\
		\bottomrule
	\end{tabular}
\end{table}

We run our adaptive algorithm for a relative tolerance of \PRreview{$\epsilon_\mathrm{rel}=1\cdot10^{-3}$}. Note that Algorithm~\ref{PR_alg:AdaptiveMIMC} is formulated in terms of an \emph{absolute} tolerance $\epsilon$. We adapt for (estimated) relative tolerances by using the current estimate for the expected value of the quantity of interest \PRreview{as a scaling factor}. The mean value of the quantity of interest was computed by our algorithm as \PRreview{$Q_L(g)=133.71$} with a standard error of \PRreview{$10.45$} in \PRreview{$L=40$} iterations. The standard deviation of the estimator is \PRreview{$0.0975$}, and the estimated bias is \PRreview{$0.0782$}, giving a total ({root mean square}) {error} (RMSE) estimate of  \PRreview{$0.125<\epsilon_\mathrm{rel}\cdot Q_L(g)$}. Figure~\ref{PR_fig:IndexSets} shows the shape of the index set for some selected iterations. We see that the adaptive algorithm mainly exploits the spatial resolutions, until the addition of more spatial levels is estimated to be too expensive \PRreview{($L=8$). After that, the approximations for the conductor and insulator material are improved up to 256 and 8\,192 terms respectively. From $L=25$ and beyond, the \emph{mixed directions} that improve the approximation for both conductor and insulator material, and the approximation for the conductor material and the mesh refinement, are activated. Observe that the final shape of the index set ($L=40$) is far from trivial, and is also not immediately representable by an anisotropic simplex.}

Finally, we investigate the performance of our adaptive method compared to standard MIMC with the common choice of \PRreview{simplices} $T_{(1,1,1)}(L)$ \PRreview{as} index sets. The mean value, RMSE and runtime for all tolerances are shown in~Table~\ref{PR_tab:RunTimeComparison}. All simulations are performed on a 2.6GHz Intel Xeon CPU with 64GB of RAM. Observe that both methods converge to the same value. The adaptive algorithm outperforms the nonadaptive MIMC method for all values of $\epsilon_\text{rel}$ considered. Note that we are not able to solve for smaller tolerances using the nonadaptive MIMC method, because of the increasing memory requirement of the available spatial resolutions.

\begin{figure}[t]
	\centering
	\includegraphics[width=0.6\textwidth]{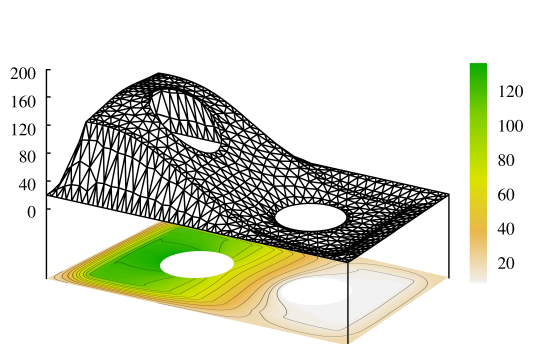}
	\caption{Mean temperature field of the heat exchanger on a 890-point mesh as an example of a nonscalar quantity of interest.}
	\label{PR_fig:MeanField}
\end{figure}

The adaptive algorithm is not limited to scalar quantities of interest. It is also possible to include multiple quantities of interest in a single simulation. We then take the worst value of the profit over all quantities considered to compute the next iterate, see~\cite{giles2015multilevel}. As an example,~Figure~\ref{PR_fig:MeanField} shows the mean value of the temperature in the heat exchanger on a mesh with 1\,524 elements, for a relative tolerance of $1\cdot10^{-3}$. The highest expected temperature is located at a point on the boundary of the hot fluid pipe, opposite to the pipe containing the coolant fluid. This is what might have been anticipated from physical considerations, assuming that the heat flux $\Phi_h$ is large enough to heat the material around the left-hand pipe to a temperature higher than $T_e$. The effect of the insulator material is obvious from the large temperature gradient present at the left side of the insulator.

\section{Discussion and Future Work}\label{PR_sec:DiscussionAndFutureWork}

We have presented a dimension-adaptive Multi-Index Monte Carlo (MIMC) method for the approximation of the expected value of a quantity of interest that is a function of the solution of a PDE with random coefficients. The method, which can be seen as a generalization of the classical MIMC method\PRreview{,} automatically finds important directions in the problem. These directions are not limited to spatial dimensions only, as is demonstrated by a numerical experiment. We have demonstrated an efficient implementation of the method,  based on a similar construction used in dimension-adaptive integration with sparse grids. 

The adaptive algorithm is particularly interesting when the \PRreview{optimal} shape of the MIMC index set is unknown or nontrivial, since it does not require a priori knowledge of the structure of the problem. In these situations, the method may include or exclude certain indices to achieve an estimator that minimizes computational effort needed to obtain a certain tolerance.

\PRreview{Finally, adaptivity can be used in combination with other techniques, such as \emph{Quasi-Monte Carlo}, see \cite{dick2013quasi} or~\cite{robbe2017multi} for the multi-index setting. We expect similar gains as outlined in this paper.}


%
\bibliographystyle{spmpsci}
%

\end{document}